\documentclass[reqno]{amsart}
\usepackage{amsmath}
\usepackage{amsthm}
\usepackage{amsfonts}
\usepackage{amssymb}
\usepackage{latexsym}
\usepackage{amscd}
\usepackage{stmaryrd}
\usepackage{color}
\usepackage{hyperref}
\usepackage{extarrows}

\usepackage{bm}

\theoremstyle{plain}
\newtheorem{thm}{Theorem}[section]
\newtheorem{cor}[thm]{Corollary}
\newtheorem{lem}[thm]{Lemma}

\numberwithin{equation}{section}

\newfont{\scyr}{wncyr10 scaled 550}

\allowdisplaybreaks
\begin{document}

\title{Topological properties of $q$-analogues of multiple zeta values}

\date{\today\thanks{The first author is supported by the National Natural Science Foundation of
China (Grant No. 11471245). The authors would like to thank the anonymous referee for the helpful comments and suggestions, which greatly improve the quality of the paper.} }

\author{Zhonghua Li \quad and \quad Ende Pan}

\address{School of Mathematical Sciences, Tongji University, No. 1239 Siping Road,
Shanghai 200092, China}

\email{zhonghua\_li@tongji.edu.cn}

\address{School of Mathematical Sciences, Tongji University, No. 1239 Siping Road,
Shanghai 200092, China}

\email{13162658945@163.com}



\begin{abstract}
In the space of bounded real-valued functions on the interval $(0,1)$, we study the convergent sequences of $q$-analogues of multiple zeta values which do not converge to $0$. And we obtain the derived sets of the set of some $q$-analogue of multiple zeta values.
\end{abstract}

\maketitle

%%---------------------------------------------------------------------------
%%------------------------Content-------------------------------------------
%%----------------------------------------------------------------------------

\section{Introduction and statement of main results}\label{Sec:Intro}

Multiple zeta values are natural generalizations of  the Riemann zeta values. Let $\mathbb{N}$ be the set of positive integers. For any $d\in\mathbb{N}$ and any multi-index $\mathbf{k}=(k_1,\ldots,k_d)\in\mathbb{N}^d$ with $k_1\geqslant 2$, the multiple zeta value $\zeta(\mathbf{k})$ is defined by the following infinite series
$$\zeta(\mathbf{k})=\zeta(k_1,\ldots,k_d)=\sum\limits_{m_1>\cdots>m_d>0}\frac{1}{m_1^{k_1}\cdots m_d^{k_d}}.$$
The condition $k_1\geqslant 2$ ensures the convergence of the above infinite series. And we call such a multi-index $\mathbf{k}$ admissible. The quantities $k_1+\cdots+k_d$ and $d$ are called weight and depth of $\mathbf{k}$, respectively. Different from other researchers' work on multiple zeta values, Kumar studied the order structure and the topological properties of the set $\mathcal{Z}$ of all multiple zeta values  in \cite{Senthil Kumar}. Taking the usual order and the usual topology of the set $\mathbb{R}$ of real numbers, Kumar computed the derived sets of the topological subspace $\mathcal{Z}$ of $\mathbb{R}$, and showed that the set $\mathcal{Z}$, ordered by $\geqslant$, is well-ordered with the order type $\omega^3$, where $\omega$ is the smallest infinite ordinal.

In this paper, we study the topological properties of some $q$-analogues of multiple zeta values. Let $q\in\mathbb{R}$ with $0<q<1$. For any $m\in\mathbb{N}$, let $[m:q]$ denote the $q$-integer
$$[m:q]=\frac{1-q^m}{1-q}=1+q+\cdots+q^{m-1}.$$
Then for any admissible multi-index $\mathbf{k}=(k_1,\ldots,k_d)\in\mathbb{N}^d$, we define the multiple $q$-zeta value $\zeta[\mathbf{k}:q]$ by
\begin{align}
&\zeta[\mathbf{k}:q]=\zeta[k_1,\ldots,k_d:q]=\sum\limits_{m_1>\cdots>m_d>0}\frac{q^{m_{1}(k_{1}-1)+\cdots+m_d(k_d-1)}}{[m_{1}:q]^{k_{1}}\cdots [m_d:q]^{k_d}}.
\label{Eq:qMZV}
\end{align}
This $q$-analogue was first studied by Bradley \cite{Bradley} and independently by Zhao \cite{Zhao1}.  Here we introduce another $q$-analogue of multiple zeta values. Let $r\in\mathbb{N}$, then we define
\begin{align}
&\zeta[\mathbf{k};r:q)=\zeta[k_1,\ldots,k_d;r:q)=\sum\limits_{m_1>\cdots>m_d>m_{d+1}>0}\frac{q^{m_{1}(k_{1}-1)+\cdots+m_d(k_d-1)}}{[m_{1}:q]^{k_{1}}\cdots [m_d:q]^{k_d}m_{d+1}^{r}}.
\label{Eq:qMZV-r}
\end{align}

Different from multiple zeta values, the multiple $q$-zeta values have a parameter $q$. Hence we work in the function space $\mathbf{B}(0,1)$, which is the set of bounded real-valued functions on the open interval $(0,1)$. Since the multiple $q$-zeta values we consider here belong to $\mathbf{B}(0,1)$ (see Corollary \ref{Cor:Belong-B}), we just study the following two subspaces of $\mathbf{B}(0,1)$:
\begin{align*}
\mathcal{QZ}=&\{\zeta[\mathbf{k}:q]\mid \mathbf{k}\text{\;is admissible}\},\\
\mathcal{QZZ}=&\{\zeta[\mathbf{k};r:q)\mid \mathbf{k}\text{\;is admissible}, r\in\mathbb{N}\}.
\end{align*}

Let $F$ be the set of functions from $(0,1)$ to $\mathbb{R}$. We define a partial order relation $\preccurlyeq$ on $F$ as follows. Let $f,g\in F$. The function $f$ is smaller than $g$,  if $f(q)<g(q)$ for any $q\in (0,1)$. We denote this by $f\prec g$. Then $f\preccurlyeq g$ if $f\prec g$ or $f=g$. We can find the maximum element of $\mathcal{QZ}$.

\begin{thm}\label{Thm:Maximum-Element}
For any admissible multi-index $\mathbf{k}$, we have $\zeta[\mathbf{k}:q]\preccurlyeq \zeta[2:q]$. In other words, $\zeta[2:q]$ is the maximum element of $\mathcal{QZ}$.
\end{thm}

While for the subspace $\mathcal{QZZ}$, we only obtain an upper bound.

\begin{thm}\label{Thm:Upper-Bound}
For any admissible multi-index $\mathbf{k}$ and any $r\in\mathbb{N}$, we have $\zeta[\mathbf{k};r:q)\prec\zeta[2:q]$. In other words, $\zeta[2:q]$ is an upper bound of $\mathcal{QZZ}$.
\end{thm}

We prove Theorem \ref{Thm:Maximum-Element} and Theorem \ref{Thm:Upper-Bound} in Section \ref{Sec:Proof-Thm1-2}.

As in \cite{Senthil Kumar}, we want to compute the derived sets of the subspace $\mathcal{QZZ}$. Hence some topology of $\mathbf{B}(0,1)$ is needed. In fact, $\mathbf{B}(0,1)$ is a complete normed space with the norm given by
$$\|f\|=\sup\limits_{q\in(0,1)}|f(q)|,\quad \text{for each\;} f\in\mathbf{B}(0,1).$$
In the following we consider $\mathbf{B}(0,1)$ with respect to the topology induced by the above norm. We are interested to find the derived sequence $(\mathcal{QZZ}^{(n)})_{n\geqslant 0}$ of the subspace $\mathcal{QZZ}$ of $\mathbf{B}(0,1)$ where $\mathcal{QZZ}^{(0)}=\mathcal{QZZ}$ and for $n\geqslant 1$, $\mathcal{QZZ}^{(n)}$ is the set of accumulation points of $\mathcal{QZZ}^{(n-1)}$ in $\mathbf{B}(0,1)$. To state our results, we further need to introduce some more notations. For an admissible multi-index $\mathbf{k}=(k_1,\ldots,k_d)\in\mathbb{N}^d$ and a nonnegative integer $n$, we set
\begin{align}
&\zeta[\mathbf{k}:q]_n=\zeta[k_1,\ldots,k_d:q]_{n}=\sum\limits_{m_1>\cdots>m_d>n}\frac{q^{m_{1}(k_{1}-1)+\cdots+m_d(k_d-1)}}{[m_{1}:q]^{k_{1}}\cdots [m_d:q]^{k_d}}.
\label{Eq:qMZV-Tails}
\end{align}
Obviously, we have $\zeta[\mathbf{k}:q]_0=\zeta[\mathbf{k}:q]$. We have the following theorem, the proof will be given in Section \ref{Sec:Proof-Thm3}.

\begin{thm}\label{Thm:Derived-Set}
We have
$$\mathcal{QZZ}^{(1)}=\{\zeta[\mathbf{k}:q]_1\mid \mathbf{k}\text{\;is admissible}\}\cup \{0\}$$
and $\mathcal{QZZ}^{(2)}=\{0\}$.
\end{thm}

Throughout the paper, the notation $\{k\}^n$ stands for $\underbrace{k,\ldots,k}_{n\text{\;terms}}$.

%%%%%%%%%%%%%%%%%%%%%%%%%%%%%%%%%%%%%%%%%%%%%%%%%%%%%%%%%%%%%%%%

\section{Proofs of  Theorem \ref{Thm:Maximum-Element} and Theorem \ref{Thm:Upper-Bound}} \label{Sec:Proof-Thm1-2}

In this section, we give proofs of Theorem \ref{Thm:Maximum-Element} and Theorem \ref{Thm:Upper-Bound}. We need the following lemmas.

\begin{lem}\label{Lem:OneTerm-Increasing}
For any $m\in\mathbb{N}$, the function $f(q)=\frac{q^m}{[m:q]}$ is monotonically increasing on the interval $(0,1)$. In particular, we have
$$0\prec\frac{q^m}{[m:q]}\prec\frac{1}{m},\quad \text{for each\;} q\in (0,1).$$
\end{lem}

\noindent {\bf Proof.}
We have
$$(1-q^m)^2f'(q)=q^{m-1}\left[m-(m+1)q+q^{m+1}\right].$$
Set $g(q)=m-(m+1)q+q^{m+1}$, then one gets
$$g'(q)=(m+1)(q^m-1)<0.$$
Hence for any $q\in (0,1)$, we have
$$g(q)>g(1)=0,$$
which induces that $f'(q)>0$.
\qed

\begin{lem}\label{Lem:Del-One-Big}
For positive integers $d,j,r,k_1,\ldots,k_d$ with $k_1\geqslant 2$ and $j\leqslant d$, and any nonnegative integer $n$, we have
\begin{align*}
&\zeta[k_1,\ldots,k_{j-1},k_{j}+1,k_{j+1},\ldots,k_d:q]\prec\zeta[k_1,\ldots,k_{j},\ldots,k_d:q],\\
&\zeta[k_1,\ldots,k_{j-1},k_{j}+1,k_{j+1},\ldots,k_d:q]_n\prec\zeta[k_1,\ldots,k_{j},\ldots,k_d:q]_n
\end{align*}
and
$$\zeta[k_1,\ldots,k_{j-1},k_{j}+1,k_{j+1},\ldots,k_d;r:q)\prec\zeta[k_1,\ldots,k_{j},\ldots,k_d;r:q).$$
We also have
$$\zeta[k_1,\ldots,k_d;r+1:q)\prec\zeta[k_1,\ldots,k_d;r:q)$$
and
$$\zeta[k_1,\ldots,k_d;r:q)\prec\zeta[k_1,\ldots,k_d,1:q].$$
\end{lem}

\noindent {\bf Proof.}
From Lemma \ref{Lem:OneTerm-Increasing}, for any $m\in\mathbb{N}$ and any $q\in (0,1)$, we have $\frac{q^m}{[m:q]}<1$. Multiplying by $\frac{q^{(k_j-1)m}}{[m:q]^{k_j}}$ on both sides, we obtain
$$\frac{q^{k_jm}}{[m:q]^{k_j+1}}<\frac{q^{(k_j-1)m}}{[m:q]^{k_j}},$$
which induces the first three inequalities stated in the lemma. For $m\in\mathbb{N}$, we have
$$\frac{1}{m^{r+1}}\leqslant \frac{1}{m^r},\quad \frac{1}{m^r}\leqslant \frac{1}{[m:q]^r}\leqslant \frac{1}{[m:q]},$$
(where the equality holds only when $m=1$) from which, the last two inequalities of the lemma follows.
\qed

\begin{lem}\label{Lem:Del-Depth-Big}
For any nonnegative integer $d$, we have $\zeta[2,\{1\}^{d+1}:q]\prec\zeta[2,\{1\}^{d}:q]$.
\end{lem}

\noindent {\bf Proof.}
We use the duality formula for multiple $q$-zeta values proved by Bradley in \cite{Bradley}: for any nonnegative integers $n$ and $m$, one has
\begin{equation}
\zeta[2+n,\{1\}^{m}:q]=\zeta[2+m,\{1\}^{n}:q].
\label{Eq:Dual-Formula}
\end{equation}
From \eqref{Eq:Dual-Formula}, to prove the lemma, it is enough to show $\zeta[d+3:q]\prec\zeta[d+2:q]$, but this follows from Lemma \ref{Lem:Del-One-Big}.
\qed

Now we prove Theorem \ref{Thm:Maximum-Element} and Theorem \ref{Thm:Upper-Bound}.

\noindent {\bf Proof of Theorem \ref{Thm:Maximum-Element}.}
Let $\mathbf{k}=(k_1,\ldots,k_d)$. By Lemmas \ref{Lem:Del-One-Big} and  \ref{Lem:Del-Depth-Big}, we have
$$\zeta[\mathbf{k}:q]=\zeta[k_1,\ldots,k_d:q]\preccurlyeq \zeta[2,\{1\}^{d-1}:q]\preccurlyeq \zeta[2:q],$$
as desired. \qed

\noindent {\bf Proof of Theorem \ref{Thm:Upper-Bound}.}
Let $\mathbf{k}=(k_1,\ldots,k_d)$. By Lemmas \ref{Lem:Del-One-Big} and \ref{Lem:Del-Depth-Big}, we have
$$\zeta[\mathbf{k};r:q)=\zeta[k_1,\ldots,k_d;r:q)\prec\zeta[2,\{1\}^{d}:q]\prec\zeta[2:q],$$
as desired. \qed

We end this section with a corollary.

\begin{cor}\label{Cor:Belong-B}
The function spaces $\mathcal{QZ}$ and $\mathcal{QZZ}$ are subsets of $\mathbf{B}(0,1)$.
\end{cor}

\noindent {\bf Proof.}
Since
$$\lim\limits_{q\rightarrow 0}\zeta[2:q]=0,\qquad \lim\limits_{q\rightarrow 1}\zeta[2:q]=\zeta(2),$$
the function $\zeta[2:q]$ is continuous on the closed interval $[0,1]$. Hence $\zeta[2:q]\in \mathbf{B}(0,1)$. Now for any admissible multi-index $\mathbf{k}$, from Theorem \ref{Thm:Maximum-Element}, we have $\zeta[\mathbf{k}:q]\in \mathbf{B}(0,1)$. Then $\mathcal{QZ}$ is a  subset of $\mathbf{B}(0,1)$. Similarly, from Theorem \ref{Thm:Upper-Bound}, we find $\mathcal{QZZ}$ is a  subset of $\mathbf{B}(0,1)$.
\qed

%%%%%%%%%%%%%%%%%%%%%%%%%%%%%%%%%%%%%%%%%%%%%%%%%%%%%%%%%%%%%%%%%%%%%%%%%%%%%%%%%%%%%%%%%%%%%%%%%%%%%%%%%%%%%%%%%%%%

\section{Proof of Theorem \ref{Thm:Derived-Set}}\label{Sec:Proof-Thm3}

In this section, we give a proof of  Theorem \ref{Thm:Derived-Set}.

\subsection{Some preliminary results}

We first compute the norms of some multiple $q$-zeta values. For this purpose, we require the following lemmas.

\begin{lem}\label{Lem:Mono-Decreasing}
For a fixed $q\in (0,1)$, the function $f(x)=\frac{xq^{x-1}}{1-q^x}$ is monotonically decreasing on the interval $[1,+\infty)$.
\end{lem}

\noindent {\bf Proof.}
We have
$$(1-q^x)^2f'(x)=q^{x-1}(1-q^x+x\log q).$$
Set $g(x)=1-q^x+x\log q$, then we have
$$g'(x)=(1-q^x)\log q<0.$$
Hence we find
$$g(x)\leqslant g(1)=1-q+\log q.$$
Set $h(x)=1-x+\log x$, where $x\in (0,1)$, then we have
$$h'(x)=-1+\frac{1}{x}=\frac{1-x}{x}>0.$$
Finally, we get
$$g(1)=h(q)<h(1)=0,$$
which implies that $f'(x)<0$.
\qed

\begin{lem}\label{Lem:Monotone}
Let $d,m_1,\ldots,m_d,k$ be positive integers such that $m_{1}\geqslant \cdots \geqslant m_{d}$ and $k \geqslant d+1$. Then the function $f(q)=\frac{q^{m_{1}(k-1)}}{{[m_{1}:q]}^{k}[m_{2}:q]\cdots [m_{d}:q]}$ is monotonically increasing on the interval $(0,1)$. In particular, for $m,k\in\mathbb{N}$ with $k\geqslant 2$, the function $\frac{q^{m(k-1)}}{{[m:q]}^{k}}$ is monotonically increasing on the interval $(0,1)$.
\end{lem}

\noindent {\bf Proof.} Since
$$f(q)=\frac{{(1-q)}^{k+d-1}q^{m_{1}(k-1)}}{(1-q^{m_{1}})^{k}(1-q^{m_{2}})\ldots(1-q^{m_{d}})},$$
taking the logarithmic derivative of $f(q)$, we get
$$\frac{f'(q)}{f(q)}=\frac{m_{1}(k-1)}{q}-\frac{k+d-1}{1-q}+\frac{k{m_{1}}{q^{m_{1}-1}}}{1-q^{m_{1}}}+\sum_{i=2}^{d}\frac{{m_{i}}{q^{m_{i}-1}}}{1-q^{m_{i}}}.$$
Using Lemma \ref{Lem:Mono-Decreasing}, we get
$$\frac{f'(q)}{f(q)}\geqslant \frac{m_{1}(k-1)}{q}-\frac{k+d-1}{1-q}+\frac{(k+d-1){m_{1}}{q^{m_{1}-1}}}{1-q^{m_{1}}},$$
which is equivalent to
$$q(1-q)(1-q^{m_1})\frac{f'(q)}{f(q)}\geqslant g(q)$$
with
$$g(q)=\left[m_{1}(k-1)(1-q)-(k+d-1)q\right](1-q^{m_1})+(k+d-1)m_{1}q^{m_{1}}(1-q).$$
Then it is enough to show that
$$g(q)>0,\quad \text{for each\;} q\in (0,1).$$

In fact, we have
\begin{align*}
&g'(q)=m_{1}-m_{1}k-k-d+1+m_{1}^{2}dq^{m_{1}-1}+(k+d-1-m_{1}d)(m_{1}+1)q^{m_{1}},\\
&g''(q)=q^{m_{1}-2}m_{1}[m_{1}(m_{1}-1)d+(m_{1}+1)(k+d-1-m_{1}d)q].
\end{align*}
Set
$$h(q)=m_{1}(m_{1}-1)d+(m_{1}+1)(k+d-1-m_{1}d)q,$$
then, since $k\geqslant d+1$, we get
$$h(0)=m_{1}(m_{1}-1)d\geqslant 0,\quad h(1)=m_{1}(k-d-1)+k+d-1>0.$$
Therefore for any $q\in (0,1)$, we have $h(q)>0$ and then $g''(q)>0$. This implies that
$$g'(q)<g'(1)=0,\quad \text{for each\;} q\in (0,1),$$
and then
$$g(q)>g(1)=0,\quad \text{for each\;} q\in (0,1).$$
This completes the proof.\qed

From Lemma \ref{Lem:Monotone},  we get the norms of height one multiple $q$-zeta values.

\begin{cor}\label{Cor:Norm-qMZV}
For any nonnegative integers $n$ and $m$, we have
$$\|\zeta[2+n,\{1\}^{m}:q]\|=\zeta(2+n,\{1\}^{m}).$$
\end{cor}

\noindent {\bf Proof}.
If $n\geqslant m$, we get the result from Lemma \ref{Lem:Monotone}. If $n\leqslant m$, applying the duality formula \eqref{Eq:Dual-Formula} and  its multiple zeta values' version, we get the result from Lemma \ref{Lem:Monotone}.
\qed

Now we provide upper and lower bounds for tails of multiple $q$-zeta values. We set
$$\Omega_1(\mathbf{k},n:q)={\left({\frac{q-1}{\log q}}\right)}^{d}{\left(\frac{q^{n+d}}{[n+d:q]}\right)}^{{k_{1}+\cdots+k_{d}-d}}\prod_{i=1}^{d}\frac{1}{k_{1}+\cdots+k_{i}-i}$$
and
$$\Omega_2(\mathbf{k},n:q)={\left({\frac{q-1}{\log q}}\right)}^{d}{\left(\frac{q^{n}}{[n:q]}\right)}^{{k_{1}+\cdots+k_{d}-d}}\prod_{i=1}^{d}\frac{1}{k_{1}+\cdots+k_{i}-i},$$
where $\mathbf{k}=(k_1,\ldots,k_d)\in\mathbb{N}^d$ and $n\in\mathbb{N}$.

\begin{lem}\label{Lem:Equivalence}
For any admissible multi-index $\mathbf{k}=(k_1,\ldots,k_d)\in\mathbb{N}^d$ and any $n\in\mathbb{N}$, we have
\begin{align}
\Omega_1(\mathbf{k},n:q)\prec\zeta[\mathbf{k}:q]_{n}\prec\Omega_2(\mathbf{k},n:q).
\label{Eq:Ineq-qMZV-Tail}
\end{align}
\end{lem}

\noindent {\bf Proof.}
We prove by induction  on  $d$. Let $q\in (0,1)$ be fixed. For $d=1$, we have to show that
\begin{align}
\frac{q-1}{\log q}\left(\frac{q^{n+1}}{[n+1:q]}\right)^{k_1-1}\frac{1}{k_1-1}<\zeta[k_1:q]_n<\frac{q-1}{\log q}\left(\frac{q^{n}}{[n:q]}\right)^{k_1-1}\frac{1}{k_1-1}.
\label{Eq:Ineq-qMZV-Tail-D1}
\end{align}
In fact, set $f_{k_1}(x)=\frac{q^{x(k_{1}-1)}}{{(1-q^{x})}^{k_{1}}}$, then we have
$$(1-q^x)^{k_1+1}f_{k_1}'(x)=q^{x(k_1-1)} (k_1-1+q^x)\log q<0,\quad (x\geqslant 1).$$
Hence $f_{k_1}(x)$ is monotonically decreasing on the interval $[1,+\infty)$ for any $k_{1}\in\mathbb{N}$. Then we obtain
$${(1-q)}^{k_{1}}\int _{n+1}^{\infty}\frac{q^{x(k_{1}-1)}}{{(1-q^{x})}^{k_{1}}}dx < \zeta[k_1:q]_n<{(1-q)}^{k_{1}}\int _{n}^{\infty}\frac{q^{x(k_{1}-1)}}{{(1-q^{x})}^{k_{1}}}dx.$$
For $k_1\geqslant 2$, the substitution $y=q^x$ gives the following identity
$${(1-q)}^{k_{1}}\int_{a}^{b}\frac{q^{x(k_{1}-1)}}{{(1-q^{x})}^{k_{1}}}dx
=\left[\frac{{(1-q)}^{k_{1}}}{(k_{1}-1)\log q}{\left(\frac{y}{1-y}\right)}^{k_{1}-1}\right]_{y=q^{a}}^{q^{b}},
$$
from which we get \eqref{Eq:Ineq-qMZV-Tail-D1}.

If $d>1$, we have
$$\zeta[\mathbf{k}:q]_n=\zeta[k_1,\ldots,k_d:q]_{n}=\sum_{m_{d}>n}\frac{q^{m_{d}(k_{d}-1)}}{{[m_{d}:q]}^{k_{d}}}\zeta[k_1,\ldots,k_{d-1}:q]_{m_{d}}.$$
Using the induction hypothesis, we get
$$\Lambda_l(q)<\zeta[\mathbf{k}:q]_n<\Lambda_r(q),$$
where
\begin{align*}
\Lambda_l(q)=&\sum_{m_{d}>n}\frac{q^{m_{d}(k_{d}-1)}}{{[m_{d}:q]}^{k_{d}}}{\left({\frac{q-1}{\log q}}\right)}^{d-1}\left(\frac{q^{m_{d}+d-1}}{{[m_{d}+d-1:q]}}\right)^{{k_{1}+\cdots+k_{d-1}-d+1}}\\
&\quad\times\prod_{i=1}^{d-1}\frac{1}{k_{1}+\cdots+k_{i}-i},
\end{align*}
and
\begin{align*}
\Lambda_r(q)=&\sum_{m_{d}>n}\frac{q^{m_{d}(k_{d}-1)}}{{[m_{d}:q]}^{k_{d}}}{\left({\frac{q-1}
{\log q}}\right)}^{d-1}\left(\frac{q^{m_{d}}}{{[m_{d}:q]}}\right)^{{k_{1}+\cdots+k_{d-1}-d+1}}\\
&\quad\times\prod_{i=1}^{d-1}\frac{1}{k_{1}+\cdots+k_{i}-i}.
\end{align*}
Since
$$\frac{q^{m_d(k_d-1)}}{[m_d:q]^{k_d}}=(1-q)^{k_d}f_{k_d}(m_d)>(1-q)^{k_d}f_{k_d}(m_d+d-1)=\frac{q^{(m_d+d-1)(k_d-1)}}{[m_d+d-1:q]^{k_d}},$$
we find
$$\Lambda_l(q)>\left(\frac{q-1}{\log q}\right)^{d-1}\prod_{i=1}^{d-1}\frac{1}{k_{1}+\cdots+k_{i}-i}\zeta[k_1+\cdots+k_d-d+1:q]_{n+d-1}.$$
Using the lower bound in the case of $d=1$, we have
$$\Lambda_l(q)>{\left({\frac{q-1}{\log q}}\right)}^{d}{\left(\frac{q^{n+d}}{[n+d:q]}\right)}^{{k_{1}+\cdots+k_{d}-d}}\prod_{i=1}^{d}\frac{1}{k_{1}+\cdots+k_{i}-i},$$
as desired. Similarly, since
$$\Lambda_r(q)=\left(\frac{q-1}{\log q}\right)^{d-1}\prod_{i=1}^{d-1}\frac{1}{k_{1}+\cdots+k_{i}-i}\zeta[k_1+\cdots+k_d-d+1:q]_{n},$$
using the upper bound in the case of $d=1$, we have
$$\Lambda_r(q)<{\left({\frac{q-1}{\log q}}\right)}^{d}{\left(\frac{q^{n}}{[n:q]}\right)}^{{k_{1}+\cdots+k_{d}-d}}\prod_{i=1}^{d}\frac{1}{k_{1}+\cdots+k_{i}-i},$$
as desired.
\qed

\subsection{Convergent sequences in $\mathcal{QZZ}$}

To prove Theorem \ref{Thm:Derived-Set}, we have to know the behaviour of the convergent sequences in the space $\mathcal{QZZ}$. We first introduce some notation. Let $(\mathbf{k}(n))_{n\in\mathbb{N}}=((k_1(n),\ldots,k_{d(n)}(n)))_{n\in\mathbb{N}}$ be a fixed sequence of admissible multi-indices. Set
$$N_2=\{n\in\mathbb{N}\mid k_1(n)=2\}\subset\mathbb{N}.$$
For any $n\in N_2$, we define
$$l(n)=\begin{cases}
i & \text{if\;} d(n)\geqslant 2 \text{\;and\;}k_2(n)=\cdots=k_{i-1}(n)=1,k_i(n)\geqslant 2,\\
1 & \text{otherwise},
\end{cases}$$
and $v(n)=d(n)-l(n)+1$. Then for $n\in N_2$ and $l(n)\geqslant 2$, we set
$$\mathbf{k}(n)=(2,\{1\}^{l(n)-2},k_{l(n)}(n),\ldots,k_{d(n)}(n))=(2,\{1\}^{l(n)-2},s_1(n),\ldots,s_{v(n)}(n)).$$
Finally, we define some subsets of $\mathbb{N}$ as follows. Let
\begin{align*}
D=&\{d(n)\mid n\in N_2\},\\
V=&\{v(n)\mid n\in N_2\},\\
W=&\{k_1(n)+\cdots+k_{d(n)}(n)\mid n\in N_2\},\\
W'=&\{s_1(n)+\cdots+s_{v(n)}(n)\mid n\in N_2,l(n)\geqslant 2\}.
\end{align*}
Then we have the following theorem.

\begin{thm}\label{Thm:Classification}
Let $(\mathbf{k}(n))_{n\in\mathbb{N}}=((k_1(n),\ldots,k_{d(n)}(n)))_{n\in\mathbb{N}}$ be a sequence of admissible multi-indices and $(r(n))_{n\in\mathbb{N}}$ be a sequence of positive integers. Assume that $0$ is not an accumulation point of the sequence $(\zeta[\mathbf{k}(n);r(n):q))_{n\in\mathbb{N}}$ in $\mathbf{B}(0,1)$. The following holds.
\begin{enumerate}
  \item [(i)] If $N_2$ is a finite set, then both the sets $\{d(n)\mid n\in\mathbb{N}\}$ and $\{k_1(n)+\cdots+k_{d(n)}(n)\mid n\in\mathbb{N}\}$ are bounded.
  \item [(ii)] Assume that $N_2$ is infinite. Then one of the following subcases holds.
  \begin{enumerate}
    \item [(ii-1)] If $D$ is bounded, then $W$ is bounded;
    \item [(ii-2)] If both $D$ and $V$ are unbounded, then there are only finitely many $n\in N_2$, such that $l(n)\geqslant 2$;
    \item [(ii-3)] If $D$ is unbounded while $V$ is bounded, then $W'$ is bounded.
  \end{enumerate}
\end{enumerate}
\end{thm}

\noindent {\bf Proof.}
(i) Assume that $N_2$ is a finite set. Without loss of generality, we may assume that $N_2=\emptyset$. Then for any $n\in\mathbb{N}$, we have $k_1(n)\geqslant 3$. Using Lemma \ref{Lem:Del-One-Big} and the duality formula \eqref{Eq:Dual-Formula}, we have
\begin{align*}
\zeta[k_{1}(n),\cdots,k_{d(n)}(n);r(n):q)\prec&\zeta[k_{1}(n),\cdots,k_{d(n)}(n),1:q]\\
\preccurlyeq&\zeta[3,\{1\}^{d(n)}:q]=\zeta[2+d(n),1:q].
\end{align*}
Taking norms, we have
$$\|\zeta[k_{1}(n),\cdots,k_{d(n)}(n);r(n):q)\|\leqslant \|\zeta[2+d(n),1:q]\|=\zeta(2+d(n),1),$$
where the last equality is from Corollary \ref{Cor:Norm-qMZV}. If $d(n)$ is unbounded, then there exists an infinite sequence $(n_l)$ of elements of $\mathbb{N}$, such that
$$\lim\limits_{l\rightarrow \infty}d(n_l)=\infty.$$
Since
$$\lim\limits_{d(n)\rightarrow \infty}\zeta(2+d(n),1)=0,$$
we get
$$\lim\limits_{l\rightarrow \infty}\|\zeta[k_{1}(n_l),\cdots,k_{d(n_l)}(n_l);r(n_l):q)\|=0,$$
which implies that $0$ is an accumulation point of the sequence $(\zeta[\mathbf{k}(n);r(n):q))_{n\in\mathbb{N}}$, a contradiction. Hence $d(n)$ is bounded.

Let $d$ be the maximal element of $\{d(n)\mid n\in\mathbb{N}\}$. For any $1\leqslant p\leqslant d$, set
$$M_p=\{n\in\mathbb{N}\mid d(n)=p\}.$$
For a fixed $p$, we show that for any $1\leqslant j\leqslant p$, the sets
$$\{k_j(n)\mid n\in M_p\}$$
are all bounded. If $M_p$ is a finite set, then the result follows easily. Now assume that $M_p$ is an infinite set.

For $j=1$, as above we have
\begin{align*}
\|\zeta[k_{1}(n),\cdots,k_p(n);r(n):q)\|\leqslant&\|\zeta[k_{1}(n),\ldots,k_p(n),1:q]\|\\
\leqslant&\|\zeta[k_{1}(n),\{1\}^{p}:q]\|=\zeta(k_1(n),\{1\}^p).
\end{align*}
If $k_{1}(n)$ is unbound for $n\in M_p$, then there exists an infinite sequence $(n_l)$ of elements of $M_p$ such that
$$\lim\limits_{l\rightarrow \infty}k_1(n_l)=\infty.$$
Therefore we have
$$\lim\limits_{l\rightarrow\infty}\|\zeta[k_{1}(n_l),\cdots,k_p(n_l);r(n_l):q)\|=0,$$
a contradiction. Hence $k_1(n)$ is bounded for $n\in M_p$.

Assume $2\leqslant j\leqslant p$. We may assume that for any $n\in M_p$, $k_j(n)\geqslant 2$. We have
\begin{align*}
&\zeta[k_{1}(n),\ldots,k_p(n);r(n):q)\\
=&\sum\limits_{m_1>\cdots>m_{j-1}>m_{j}}\prod_{i=1}^{j-1}\frac{q^{m_{i}(k_{i}(n)-1)}}{[m_{i}:q]^{k_{i}(n)}}
\sum\limits_{m_{j}>\cdots>m_p>m_{p+1}>0}\prod_{i=j}^{p}\frac{q^{m_{i}(k_{i}(n)-1)}}{[m_{i}:q]^{k_{i}(n)}}\frac{1}{m_{p+1}^{r(n)}}\\
\prec&\zeta[k_{1}(n),\ldots,k_{j-1}(n):q]\zeta[k_{j}(n),\ldots,k_p(n);r(n):q).
\end{align*}
Assume that we have shown $k_1(n),\ldots,k_{j-1}(n)$ are all bounded for $n\in M_p$, then $\{\|\zeta[k_{1}(n),\ldots,k_{j-1}(n):q]\|\mid n\in M_p\}$ is bounded. If $k_j(n)$ is unbounded for $n\in M_p$, then as in the case of $j=1$, there exists an infinite sequence $(n_l)$ of elements of $M_p$, such that
$$\lim\limits_{l\rightarrow\infty}\|\zeta[k_{j}(n_l),\ldots,k_p(n_l);r(n_l):q)\|=0.$$
Therefore, we again get a contradiction. Hence $k_j(n)$ is bounded for $n\in M_p$.

Finally, we find that the weight of $\mathbf{k}(n)$'s are bounded, and hence the proof of (i) follows.

\noindent Proof of (ii). Assume that $N_2$ is an infinite set.

\noindent Proof of (ii-1). If $D$ is bounded, set $d=\max D$. We only need to prove that for any $2\leqslant j\leqslant d$, $k_j(n)$ is bounded for $n\in N_2$. Then one may use a similar argument as in (i) to get the result.

\noindent Proof of (ii-2). Assume that both $D$ and $V$ are unbounded and there are infinitely many $n\in N_2$, such that $l(n)\geqslant 2$. Without loss of generality, we may assume that for any $n\in N_2$, $l(n)\geqslant 2$.  Then for $n\in N_2$, we have
\begin{align*}
&\zeta[\mathbf{k}(n);r(n):q)=\zeta[k_{1}(n),\ldots,k_{l(n)+v(n)-1}(n);r(n):q)\\
=&\sum\limits_{m_1>\cdots>m_{l(n)-1}>m_{l(n)}}\prod_{i=1}^{l(n)-1}\frac{q^{m_{i}(k_{i}(n)-1)}}{[m_{i}:q]^{k_{i}(n)}}\\
&\times\sum\limits_{m_{l(n)}>\cdots>m_{l(n)+v(n)-1}>m_{l(n)+v(n)}>0}\prod_{i=l(n)}^{l(n)+v(n)-1}\frac{q^{m_{i}(k_{i}(n)-1)}}{[m_{i}:q]^{k_{i}(n)}}\frac{1}{m_{l(n)+v(n)}^{r(n)}}.
\end{align*}
For $m_{l(n)}>\cdots>m_{l(n)+v(n)-1}>m_{l(n)+v(n)}>0$, we have $m_{l(n)}\geqslant v(n)+1>v(n)$. Hence
\begin{align*}
&\zeta[\mathbf{k}(n);r(n):q)\prec\sum\limits_{m_1>\cdots>m_{l(n)-1}>v(n)}\prod_{i=1}^{l(n)-1}\frac{q^{m_{i}(k_{i}(n)-1)}}{[m_{i}:q]^{k_{i}(n)}}\\
&\times\sum\limits_{m_{l(n)}>\cdots>m_{l(n)+v(n)-1}>m_{l(n)+v(n)}>0}\prod_{i=l(n)}^{l(n)+v(n)-1}\frac{q^{m_{i}(k_{i}(n)-1)}}{[m_{i}:q]^{k_{i}(n)}}\frac{1}{m_{l(n)+v(n)}^{r(n)}}\\
=&\zeta[k_1(n),\ldots,k_{l(n)-1}(n):q]_{v(n)}\zeta[k_{l(n)}(n),\ldots,k_{l(n)+v(n)-1}(n);r(n):q).
\end{align*}
Using Lemma \ref{Lem:Del-One-Big} and the duality formula \eqref{Eq:Dual-Formula}, we have
\begin{align*}
\zeta[\mathbf{k}(n);r(n):q)\prec&\zeta[2,\{1\}^{l(n)-2}:q]_{v(n)}\zeta[2,\{1\}^{v(n)}:q]\\
=&\zeta[2,\{1\}^{l(n)-2}:q]_{v(n)}\zeta[v(n)+2:q].
\end{align*}
By Lemma \ref{Lem:Equivalence}, we get
$$\zeta[\mathbf{k}(n);r(n):q)\prec\left(\frac{q-1}{\log q}\right)^{l(n)-1}\frac{q^{v(n)}}{[v(n):q]}\zeta[v(n)+2:q].$$
Using Lemma \ref{Lem:OneTerm-Increasing}, Corollary \ref{Cor:Norm-qMZV} and the fact that the function $\frac{q-1}{\log q}$ is monotonically increasing on $(0,1)$, we have
$$\|\zeta[\mathbf{k}(n);r(n):q)\|\leqslant \frac{1}{v(n)}\zeta(v(n)+2).$$
Then from the unboundedness of $V$, we get a contradiction.

\noindent Proof of (ii-3). Assume that $D$ is unbounded and $V$ is bounded. Then $l(n)$ is unbounded for $n\in N_2$. Hence
$$\widetilde{N_2}=\{n\in N_2\mid l(n)\geqslant 2\}$$
is an infinite subset of $N_2$. Set $v=\max V$, and for $1\leqslant p\leqslant v$, set
$$\widetilde{M_p}=\{n\in \widetilde{N_2}\mid v(n)=p\}.$$
For a fixed $p$, we need to show that $s_1(n),\ldots,s_p(n)$ are all bounded for $n\in \widetilde{M_p}$. Now
\begin{align*}
\zeta[\mathbf{k}(n);r(n):q)\prec&\zeta[2,\{1\}^{l(n)-2}:q]\zeta[s_1(n),\ldots,s_p(n);r(n):q)\\
=&\zeta[l(n):q]\zeta[s_1(n),\ldots,s_p(n);r(n):q).
\end{align*}
Then similarly as in the proof of (i), we get the result.
\qed

As a consequence, we get the following result, which is used to compute $\mathcal{QZZ}^{(1)}$.

\begin{cor}\label{Cor:Classification}
Let $(\mathbf{k}(n))_{n\in\mathbb{N}}=((k_{1}(n),\ldots,k_{d(n)}(n)))_{n\in\mathbb{N}}$ be a sequence of admissible multi-indices and $(r(n))_{n\in\mathbb{N}}$ be a sequence of positive integers. Assume that for any $n_1\neq n_2$, $\zeta[\mathbf{k}(n_1);r(n_1):q)\neq \zeta[\mathbf{k}(n_2);r(n_2):q)$. If $0$ is not an accumulation point of the sequence $(\zeta[\mathbf{k}(n);r(n):q))_{n\in\mathbb{N}}$ in $\mathbf{B}(0,1)$, then $((\mathbf{k}(n),r(n)))_{n\in \mathbb{N}}$ has a subsequence of one of the following types:
\begin{align}
&((k_{1},\ldots,k_{d},\varphi(n)+2))_{n\in \mathbb{N}},
\label{Eq:Type-1}\\
&((2,\{1\}^{\psi(n)},r))_{n\in \mathbb{N}},
\label{Eq:Type-2}\\
&((2,\{1\}^{\psi(n)},\varphi(n)+2))_{n\in \mathbb{N}},
\label{Eq:Type-3}\\
&((2,\{1\}^{\psi(n)},k_{1},\ldots,k_{d},r))_{n\in \mathbb{N}},
\label{Eq:Type-4}\\
&((2,\{1\}^{\psi(n)},k_{1},\ldots,k_{d},\varphi(n)+2))_{n\in \mathbb{N}},
\label{Eq:Type-5}
\end{align}
where $(k_{1},\ldots,k_{d})$ is a fixed admissible multi-index, $r$ is a fixed positive integer and ${(\psi(n))}_{n\in \mathbb{N}}$, ${(\varphi(n))}_{n\in \mathbb{N}}$ are strictly
increasing sequences in $\mathbb{N}$.
\end{cor}

\noindent {\bf Proof.}
We use the same notation as in Theorem \ref{Thm:Classification}. If $N_2$ is finite, then by Theorem \ref{Thm:Classification}, $d(n)$ and $k_1(n)+\cdots+k_{d(n)}(n)$ are bounded for $n\in\mathbb{N}$. Hence there exists an infinite subset $A$ of $\mathbb{N}$, such that $d(n)=d$ is a constant for any $n\in A$. Since $k_1(n)$ is bounded for $n\in A$, there exists an infinite subset $A_1$ of $A$, such that $k_1(n)=k_1$ is a constant for any $n\in A_1$. Similarly, there exists an infinite subset $A_2$ of $A_1$, such that $k_2(n)=k_2$ is a constant for any $n\in A_2$. And finally, there exists an infinite subset $B$ of $A$, such that
$$k_1(n)=k_1,\ldots,k_d(n)=k_d$$
are all constants for any $n\in B$. Now $(r(n))_{n\in B}$ must be unbounded, hence $((\mathbf{k}(n),r(n)))_{n\in \mathbb{N}}$ has a subsequence of the form \eqref{Eq:Type-1}.

Now assume that $N_2$ is an infinite set. If $D$ is bounded, then by Theorem \ref{Thm:Classification}, $W$ is bounded. A similar argument as above implies that $((\mathbf{k}(n),r(n)))_{n\in \mathbb{N}}$ has a subsequence of the form \eqref{Eq:Type-1}. If both $D$ and $V$ are unbounded, then by Theorem \ref{Thm:Classification}, there is an infinite subset $A$ of $N_2$, such that $l(n)=1$ for all $n\in A$. Then $((\mathbf{k}(n),r(n)))_{n\in \mathbb{N}}$ has a subsequence of the form \eqref{Eq:Type-2} or of the form \eqref{Eq:Type-3} according to the sequence $(r(n))_{n\in A}$ is bounded or unbounded. Finally, if $D$ is unbounded while $V$ is bounded, then by Theorem \ref{Thm:Classification}, $((\mathbf{k}(n),r(n))_{n\in \mathbb{N}})$ has a subsequence of the form \eqref{Eq:Type-4} or of the form \eqref{Eq:Type-5}.
\qed

To determine $\mathcal{QZZ}^{(2)}$, we need the following.

\begin{thm}\label{Thm:qMZVTail-Classification}
Let $(\mathbf{k}(n))_{n\in\mathbb{N}}=((k_{1}(n),\ldots,k_{d(n)}(n)))_{n\in\mathbb{N}}$ be a sequence of admissible multi-indices. Assume that for any $n_1\neq n_2$, $\zeta[\mathbf{k}(n_1):q]_1\neq \zeta[\mathbf{k}(n_2):q]_1$. If $0$ is not an accumulation point of the sequence $(\zeta[\mathbf{k}(n):q]_1)_{n\in\mathbb{N}}$ in $\mathbf{B}(0,1)$, then $(\mathbf{k}(n))_{n\in \mathbb{N}}$ has a subsequence of one of the following types:
\begin{align}
&((2,\{1\}^{\psi(n)}))_{n\in \mathbb{N}},
\label{Eq:qType-1}\\
&((2,\{1\}^{\psi(n)},k_{1},\ldots,k_{d}))_{n\in \mathbb{N}},
\label{Eq:qType-2}
\end{align}
where $(k_{1},\ldots,k_{d})$ is a fixed admissible multi-index and ${(\psi(n))}_{n\in \mathbb{N}}$ is a strictly
increasing sequence in $\mathbb{N}$.
\end{thm}

\noindent {\bf Proof.}
We can prove similarly as in Theorem \ref{Thm:Classification} and Corollary \ref{Cor:Classification}. While since for $k_{d(n)}(n)\geqslant 2$ it holds
$$\zeta[k_{1}(n),\ldots,k_{d(n)}(n):q]_{1}\prec\zeta[k_{1}(n),\ldots,k_{d(n)-1}(n):q]\zeta[k_{d(n)}(n):q]_{1},$$
if we have shown $k_1(n),\ldots,k_{d(n)-1}(n)$ are all bounded, then $k_{d(n)}(n)$ is also bounded .
\qed

\subsection{Some lemmas}

To prove Theorem \ref{Thm:Derived-Set}, we recall the concept of double tails of multiple zeta values of Akhilesh from \cite{Akhilesh}. Let $\mathbf{k}=(k_1,\ldots,k_d)\in\mathbb{N}^d$ be an admissible multi-index and $n,p$ be two nonnegative integers. Then we define
$$\zeta(\mathbf{k})_{p,n}=\zeta(k_1,\ldots,k_d)_{p,n}=\sum\limits_{m_1>\cdots>m_d>n}\binom{m_{1}+p}{p}^{-1}\frac{1}{m_1^{k_1}\cdots m_d^{k_d}}.$$
We need the duality formula of double tails of multiple zeta values. Any admissible multi-index has the form
$$\mathbf{k}=(a_1+1,\{1\}^{b_1-1},\ldots,a_s+1,\{1\}^{b_s-1}),$$
where $s,a_1,b_1,\ldots,a_s,b_s\in\mathbb{N}$. Then the dual index of $\mathbf{k}$ is defined as
$$\overline{\mathbf{k}}=(b_s+1,\{1\}^{a_s-1},\ldots,b_1+1,\{1\}^{a_1-1}).$$

\begin{lem}[\cite{Akhilesh}]\label{Lem:DoubleTail}
Let $\mathbf{k}$ be an admissible multi-index and $\overline{\mathbf{k}}$ be its dual. Then for any nonnegative integers $p$ and $n$, we have
\begin{align}
\zeta(\mathbf{k})_{p,n}=\zeta(\overline{\mathbf{k}})_{n,p}.
\label{Eq:DoubleTail-Duality}
\end{align}
\end{lem}

Let $p=n=0$ in \eqref{Eq:DoubleTail-Duality}, we get the well-known duality formula of multiple zeta values
$$\zeta(\mathbf{k})=\zeta(\overline{\mathbf{k}}).$$

We need to compare the values of the double tails of multiple zeta values.

\begin{lem}\label{Lem:Compare-DoubleTail}
For any $d,j,k_1,\ldots,k_d\in\mathbb{N}$ with $k_1\geqslant 2$ and $j\leqslant d$, we have
$$\zeta(k_1,\ldots,k_{j-1},k_{j}+1,k_{j+1},\ldots,k_d)_{0,1}<\zeta(k_1,\ldots,k_{j},\ldots,k_d)_{0,1}$$
and
$$\zeta(2,\{1\}^{d})_{0,1}<\zeta(2,\{1\}^{d-1})_{0,1}.$$
\end{lem}

\noindent {\bf Proof.}
The first inequality follows from
$$\frac{1}{m_j^{k_j+1}}<\frac{1}{m_j^{k_j}}$$
for positive integer  $m_j>1$. For the second inequality, using the duality formula \eqref{Eq:DoubleTail-Duality}, we have
\begin{align*}
\zeta(2,\{1\}^{d})_{0,1}=&\zeta(d+2)_{1,0}=\sum\limits_{m=1}^\infty\binom{m+1}{1}^{-1}\frac{1}{m^{d+2}}\\
<&\sum\limits_{m=1}^\infty\binom{m+1}{1}^{-1}\frac{1}{m^{d+1}}=\zeta(d+1)_{1,0}=\zeta(2,\{1\}^{d-1})_{0,1},
\end{align*}
as desired.
\qed

Finally, to show some sequence of $\mathbf{B}(0,1)$ does not converge, we need the following simple result.

\begin{lem}\label{lem:Uniform-Convergence}
Let the sequence ${{(f_n)}_{n\in \mathbb{N}}}$ converge to $f$ in $\mathbf{B}(0,1)$ as $n$ tends to infinity.  Then $\|f_n\|$ is convergent to $\|f\|$, and for any $q\in (0,1)$, $f_n(q)$ is convergent to $f(q)$ in $\mathbb{R}$.
\end{lem}

\noindent {\bf Proof.}
We have
$$\lim_{n\rightarrow\infty}\parallel f_n-f \parallel=0.$$
Then the facts
$$|\|f_n\|-\|f\|| \leqslant \|f_n-f \|$$
and
$$|f_n(q)-f(q)|\leqslant \|f_n-f \|,\quad (\text{for each\;} q\in (0,1))$$
imply the results.
\qed

\subsection{Proof of Theorem \ref{Thm:Derived-Set}}
Now we prove Theorem \ref{Thm:Derived-Set}. We first compute $\mathcal{QZZ}^{(1)}$. Let $n\in\mathbb{N}$. Using Lemma \ref{Lem:Del-One-Big} and the duality formula \eqref{Eq:Dual-Formula}, we have
$$\zeta[3,\{1\}^{n-1};1:q)\prec\zeta[3,\{1\}^n:q]=\zeta[n+2,1:q].$$
Then by Corollary \ref{Cor:Norm-qMZV}, we have
$$\|\zeta[3,\{1\}^{n-1};1:q)\|\leqslant \|\zeta[n+2,1:q]\|=\zeta(n+2,1).$$
Since $\lim\limits_{n\rightarrow \infty}\zeta(n+2,1)=0$, we get
$$\lim\limits_{n\rightarrow\infty}\|\zeta[3,\{1\}^{n-1};1:q)\|=0,$$
which implies that $0\in \mathcal{QZZ}^{(1)}$.

Similarly, let $\mathbf{k}=(k_1,\ldots,k_d)$ be an admissible multi-index and $n\in\mathbb{N}$. We have
\begin{align*}
&\zeta[k_{1},\ldots,k_{d};n+2:q)-\zeta[k_1,\ldots,k_d:q]_{1}\\
=&\sum_{m_{d+1}=2}^{\infty}\sum\limits_{m_1>\cdots>m_{d}>m_{d+1}}\prod_{i=1}^{d}\frac{q^{m_{i}(k_{i}-1)}}{[m_{i}:q]^{k_{i}}}\cdot
\frac{1}{m_{d+1}^{n+2}}\\
\prec&\zeta[k_1,\ldots,k_d:q]\sum\limits_{m=2}^\infty\frac{1}{m^{n+2}},
\end{align*}
which implies that
$$\|\zeta[k_{1},\ldots,k_{d};n+2:q)-\zeta[k_1,\ldots,k_d:q]_{1}\|\leqslant \|\zeta[k_1,\ldots,k_d:q]\|\sum\limits_{m=2}^\infty\frac{1}{m^{n+2}}.$$
Since
$$\lim\limits_{n\rightarrow\infty}\sum\limits_{m=2}^\infty\frac{1}{m^{n+2}}=0,$$
we get
$$\lim\limits_{n\rightarrow\infty}\zeta[k_{1},\ldots,k_{d};n+2:q)=\zeta[k_1,\ldots,k_d:q]_{1}.$$
And then $\zeta[k_1,\ldots,k_d:q]_{1}\in\mathcal{QZZ}^{(1)}$.

Conversely, for any $f\in\mathcal{QZZ}^{(1)}$, there exists a sequence $(\zeta[\mathbf{k}(n),r(n)))_{n\in\mathbb{N}}$ such that $\mathbf{k}(n)$ is admissible,  $r(n)\in\mathbb{N}$ and
$$\lim\limits_{n\rightarrow \infty}\zeta[\mathbf{k}(n);r(n):q)=f.$$
We may assume that $f\neq 0$ and for any $n_1\neq n_2$,  $\zeta[\mathbf{k}(n_1);r(n_1):q)\neq \zeta[\mathbf{k}(n_2);r(n_2):q)$. By Corollary \ref{Cor:Classification}, the sequence $((\mathbf{k}(n),r(n)))_{n\in\mathbb{N}}$ has a subsequence of one of the types \eqref{Eq:Type-1}-\eqref{Eq:Type-5}. Without loss of generality, we may assume that the sequence $((\mathbf{k}(n),r(n)))_{n\in\mathbb{N}}$ itself is one of  the types \eqref{Eq:Type-1}-\eqref{Eq:Type-5}.

Let $((\mathbf{k}(n),r(n)))_{n\in\mathbb{N}}$ be of the type \eqref{Eq:Type-1}. Similarly as above, we have
$$\lim\limits_{n\rightarrow\infty}\|\zeta [k_{1},\ldots,k_{d};\varphi(n)+2:q)-\zeta[k_1,\ldots,k_d:q]_{1}\|=0.$$
Therefore, $f=\zeta[k_1,\ldots,k_d:q]_{1}$.

Let $((\mathbf{k}(n),r(n)))_{n\in\mathbb{N}}$ be one of the types \eqref{Eq:Type-2}-\eqref{Eq:Type-5}. On the one hand, using Lemma \ref{lem:Uniform-Convergence}, we have
$$\|f\|=\lim\limits_{n\rightarrow\infty}\|\zeta[\mathbf{k}(n);r(n):q)\|.$$
As the function $\zeta[\mathbf{k}(n);r(n):q)$ is continuous in the interval $[0,1]$, by the definition of norms, we have
$$\|\zeta[\mathbf{k}(n);r(n):q)\|\geqslant \zeta(\mathbf{k}(n),r(n)).$$
Using \cite[Theorem 3]{Senthil Kumar}, for $n$ large enough, we have
$$\|\zeta[\mathbf{k}(n);r(n):q)\|\geqslant \zeta(2,\{1\}^{\psi(n)},k_1,\ldots,k_d,\varphi(n)+2).$$
Here if $((\mathbf{k}(n),r(n)))_{n\in\mathbb{N}}$ is of the type \eqref{Eq:Type-2} or \eqref{Eq:Type-3}, we take $(k_1,\ldots,k_d)$ to be any fixed admissible multi-index. And if $((\mathbf{k}(n),r(n)))_{n\in\mathbb{N}}$ is of the type \eqref{Eq:Type-2} or \eqref{Eq:Type-4}, we take $(\varphi(n))_{n\in\mathbb{N}}$ to be any fixed strictly increasing sequence in $\mathbb{N}$. Now since
\begin{align*}
&\zeta(2,\{1\}^{\psi(n)},k_1,\ldots,k_d,\varphi(n)+2)\\
=&\zeta(2,\{1\}^{\psi(n)},k_1,\ldots,k_d)_{0,1}+\zeta(2,\{1\}^{\psi(n)},k_1,\ldots,k_d,\varphi(n)+2)_{0,1}
\end{align*}
and
\begin{align*}
&\zeta(2,\{1\}^{\psi(n)},k_1,\ldots,k_d,\varphi(n)+2)_{0,1}\\
<&\zeta(2,\{1\}^{\psi(n)},k_1,\ldots,k_d)\zeta(\varphi(n)+2)_{0,1}\\
<&\zeta(2,\{1\}^{\psi(n)})\zeta(\varphi(n)+2)_{0,1}\\
=&\zeta(\psi(n)+2)\zeta(\varphi(n)+2)_{0,1}\longrightarrow 0\quad (n\rightarrow \infty),
\end{align*}
we obtain
$$\lim_{n\rightarrow \infty}\zeta(2,\{1\}^{\psi(n)},k_1,\ldots,k_d,\varphi(n)+2)=\lim\limits_{n\rightarrow\infty}\zeta(2,\{1\}^{\psi(n)},k_1,\ldots,k_d)_{0,1}.$$
Using the duality formula \eqref{Eq:DoubleTail-Duality}, we get
$$\|f\|\geqslant \lim\limits_{n\rightarrow\infty} \zeta(\overline{\mathbf{k}},\psi(n)+2)_{1,0}=\zeta(\overline{\mathbf{k}})_{1,1}>0,$$
where $\overline{\mathbf{k}}$ is the dual index of $(k_1,\ldots,k_d)$. On the other hand, by Lemmas \ref{Lem:Del-One-Big} and \ref{Lem:Del-Depth-Big}, we have
$$0\prec\zeta[\mathbf{k}(n);r(n):q)\prec\zeta[2,\{1\}^{\psi(n)+1}:q]=\zeta[\psi(n)+3:q].$$
Since for a fixed $q\in (0,1)$, it holds
\begin{align*}
\zeta[\psi(n)+3:q]=&q^{\psi(n)+2}+\sum\limits_{m=2}^\infty\frac{q^{(\psi(n)+2)m}}{[m:q]^{\psi(n)+3}}\\
<&q^{\psi(n)+2}+\sum\limits_{m=2}^\infty\frac{1}{m^{\psi(n)+3}}\rightarrow 0,\quad (n\rightarrow \infty),
\end{align*}
we find $f(q)=0$ for any $q\in (0,1)$ from Lemma \ref{lem:Uniform-Convergence}. Hence the sequence $(\zeta[\mathbf{k}(n);r(n):q))_{n\in\mathbb{N}}$ does not converge in $\mathbf{B}(0,1)$ provided  that $((\mathbf{k}(n),r(n)))_{n\in\mathbb{N}}$ is one of the types \eqref{Eq:Type-2}-\eqref{Eq:Type-5}.

Summarily, we have
$$\mathcal{QZZ}^{(1)}=\{\zeta[\mathbf{k}:q]_1\mid \mathbf{k}\text{\;is admissible}\}\cup \{0\}.$$

Now we compute $\mathcal{QZZ}^{(2)}$. Since
$$\zeta[3,\{1\}^n:q]_1\prec\zeta[3,\{1\}^n:q]=\zeta[n+2,1:q],$$
we have
$$\lim\limits_{n\rightarrow\infty}\|\zeta[3,\{1\}^n:q]_1\|=0.$$
And hence $0\in \mathcal{QZZ}^{(2)}$.

Conversely, for any $f\in\mathcal{QZZ}^{(2)}$, there exists a sequence $(\zeta[\mathbf{k}(n):q]_1)_{n\in\mathbb{N}}$ such that $\mathbf{k}(n)$ is admissible and
$$\lim\limits_{n\rightarrow \infty}\zeta[\mathbf{k}(n):q]_1=f.$$
We may assume that $f\neq 0$ and for any $n_1\neq n_2$,  $\zeta[\mathbf{k}(n_1):q]_1\neq \zeta[\mathbf{k}(n_2):q]_1$. By Theorem \ref{Thm:qMZVTail-Classification}, the sequence $(\mathbf{k}(n))_{n\in\mathbb{N}}$ has a subsequence of the type \eqref{Eq:qType-1} or of the type \eqref{Eq:qType-2}. Without loss of generality, we may assume that the sequence $(\mathbf{k}(n))_{n\in\mathbb{N}}$ itself is of the type \eqref{Eq:qType-1} or of the type \eqref{Eq:qType-2}.

One the one hand, using Lemma \ref{Lem:Compare-DoubleTail} and the duality formula \eqref{Eq:DoubleTail-Duality}, we have
$$\|\zeta[\mathbf{k}(n):q]_1\|\geqslant \zeta(\mathbf{k}(n))_{0,1}\geqslant \zeta(2,\{1\}^{\psi(n)},k_1,\ldots,k_d)_{0,1}=\zeta(\overline{\mathbf{k}},\psi(n)+2)_{1,0},$$
where $\overline{\mathbf{k}}$ is the dual index of $(k_1,\ldots,k_d)$. Here if $(\mathbf{k}(n))_{n\in\mathbb{N}}$ is of the type \eqref{Eq:qType-1}, we take $(k_1,\ldots,k_d)$ to be any fixed admissible multi-index. Therefore, we find
$$\|f\|=\lim\limits_{n\rightarrow\infty}\|\zeta[\mathbf{k}(n):q]_1\|\geqslant \lim\limits_{n\rightarrow\infty}\zeta(\overline{\mathbf{k}},\psi(n)+2)_{1,0}=\zeta(\overline{\mathbf{k}})_{1,1}>0.$$
On the other hand, using Lemmas \ref{Lem:Del-One-Big} and \ref{Lem:Del-Depth-Big}, we have
$$0\prec \zeta[\mathbf{k}(n):q]_1\prec \zeta[\mathbf{k}(n):q]\preccurlyeq \zeta[2,\{1\}^{\psi(n)}:q]=\zeta[\psi(n)+2:q].$$
While as shown above, for any fixed $q\in (0,1)$, it holds.
$$\lim\limits_{n\rightarrow \infty}\zeta[\psi(n)+2:q]=0.$$
Therefore $f(q)=0$ for any $q\in (0,1)$ from Lemma \ref{lem:Uniform-Convergence}. Hence $(\zeta[\mathbf{k}(n):q]_1)_{n\in\mathbb{N}}$ does not converge in $\mathbf{B}(0,1)$.

Finally, we get $\mathcal{QZZ}^{(2)}=\{0\}$. And Theorem \ref{Thm:Derived-Set} is proved. \qed

\end{document}